\documentclass[11pt]{article}
\usepackage[margin=2cm]{geometry}
\usepackage[T1]{fontenc}
\usepackage{natbib}
\usepackage[figuresright]{rotating}
\usepackage{float}
\usepackage{graphicx}
\usepackage{epstopdf}
\usepackage{url}
\usepackage[colorlinks,citecolor=blue,linkcolor=blue,urlcolor=red,pagebackref]{hyperref}
\usepackage{algorithmic}
\usepackage{algorithm}
\usepackage{amsmath,amsthm}
\usepackage{longtable}

\theoremstyle{definition}

\theoremstyle{remark}

\providecommand{\keywords}[1]
{
  \small	
  \textbf{\textit{Keywords---}} #1
}

\begin{document}

\title{Variable Neighborhood Search for the Multi-Depot Multiple Set Orienteering Problem}

\author{Ravi Kant \\ 
\small Department of Computer Science and Information Systems \\
\small Birla Institute of Technology and Science, Pilani \\
\small Pilani-333031, Rajasthan, India \\
\small \texttt{p20190020@pilani.bits-pilani.ac.in} \and 
Salmaan Shahid \\ 
\small Gameberry Labs \\ 
\small Bengaluru-560103, Karnataka, India \\ 
\small \texttt{salmaanshahid25@gmail.com} \and 
Anuvind Bhat \\ 
\small School of Computer Science \\
\small Carnegie Mellon University \\
\small 5000 Forbes Avenue \\
\small Pittsburgh, PA 15213, USA \\ 
\small \texttt{anuvind@cmu.edu} \and 
Abhishek Mishra \\ 
\small Department of Computer Science and Information Systems \\
\small Birla Institute of Technology and Science, Pilani \\
\small Pilani-333031, Rajasthan, India \\ 
\small \texttt{abhishek.mishra@pilani.bits-pilani.ac.in}}

\maketitle

\begin{abstract}
This paper introduces a variant of the Set Orienteering Problem (SOP), the multi-Depot multiple Set Orienteering Problem (mDmSOP). It generalizes the SOP by grouping nodes into mutually exclusive sets (clusters) with associated profits. Profit can be earned if any node within the set is visited. Multiple travelers, denoted by $t$ $(> 1)$, are employed, with each traveler linked to a specific depot. The primary objective of the problem is to maximize profit collection from the sets within a predefined budget. A novel formulation is introduced for the mDmSOP. The paper utilizes the Variable Neighborhood Search (VNS) meta-heuristic to solve the mDmSOP on small, medium, and large instances from the Generalized Traveling Salesman Problem (GTSP) benchmark. The results demonstrate the VNS's superiority in robustness and solution quality, as it requires less computational time than solving the mathematical formulation with GAMS 37.1.0 and CPLEX. Additionally, increasing the number of travelers leads to significant improvements in profits.
\end{abstract}


\maketitle

\section{Introduction}\label{sec:Introduction}
The routing problems with the profits category enclose diverse challenges connected by a specified objective function, unlike classical routing problems, where not every customer must be serviced. Instead, each customer is typically associated with a profit value, and the task becomes selecting the optimal subset of customers to serve while adhering to various constraints and optimizing a specified objective function (such as maximizing total profit, minimizing travel costs, or maximizing profit-cost differentials). One extensively studied problem over the past few decades is the Orienteering Problem (OP) introduced by \cite{golden1987orienteering}, falling into the node routing problem category. Unlike classical routing problems, the OP does not mandate visiting each node. Instead, each node is typically assigned a profit value. The objective is to select the optimal subset of nodes to serve using a single traveler, aiming to maximize profits within a specified upper bound (i.e., distance, time, etc.). A multitude of meta-heuristics and exact algorithms have been proposed, as documented in the work of \cite{golden1988multifaceted}, \cite{chao1996fast}, \cite{liang2002meta}, \cite{chekuri2012improved}, \cite{kobeaga2018efficient}, and \cite{fischetti1997branch} for the OP. The extensive survey for this is done by \cite{vansteenwegen2011orienteering} and \cite{archetti2015chapter}. After that, various variants of OP are proposed and studied by \cite{chao1996team}, \cite{angelelli2014clustered}, \cite{angelelli2017probabilistic}, \cite{archetti2018set} and other researchers. A new variant of the OP called the Set Orienteering Problem (SOP) has been proposed by \cite{archetti2018set}, where the nodes are divided into mutually exclusive clusters, and the profit associated with clusters can only be gained if at least one node is visited by the traveler. It is based on profit associated with clusters rather than nodes. Some of the solutions for the SOP are proposed by \cite{archetti2018set}, \cite{pvenivcka2019variable}, and \cite{carrabs2021biased}, which uses the combination of the Lin–Kernighan heuristic and Tabu-search meta-heuristic, the VNS and the Biased Random-Key Genetic Algorithm (BRKGA) respectively. It is observed that the BRKGA meta-heuristic performs better with an average gap of 0.19\% and with a 1.05\% deviation solution percentage gap calculated on the percentage gap between solution and solution average than the other two algorithms described above. The BRKGA is also better in terms of stability.

In this paper, we are particularly interested in suggesting a meta-heuristics for the generalized variant of the Set Orienteering Problem called the multi-Depot multiple Set Orienteering Problem (mDmSOP) proposed by \cite{icores24}, where customers (vertices) are grouped in mutually disjoint sets. The profit is associated with each cluster (set) rather than each customer. We use more than one traveler and the same number of depots associated with the travelers to collect the profit and maximize the profit within a fixed upper bound. We study the mDmSOP for cumulative budget where the traveler can share their budget as well as individual budget for each traveler.

The SOP is important to study as it has some practical applications suggested by the above researchers, but there can be many problems that cannot be modeled using only a single traveler, like waste management systems in smart cities and goods supply where more than one distributors are needed. The SOP found the application in the supply chain where customers belonging to different supply chains are grouped together, but in mDmSOP, consider the case where more than one distribution point is available, and the distributors can gain the maximum profit by providing the service to the chain optimal to reach for them, by this way the distributors may provide the low price to the customers in the supply chain by serving the complete demand of all the customers to a single customer in a group by indirectly serving all the customer's needs in a chain.

The following is a summary of our contribution to this paper. In section \ref{section:2}, we define the mDmSOP formally and give the mathematical formulation. We then propose a meta-heuristic in section \ref{section:3}, combining Greedy, Hungarian, and VNS procedures to solve the mDmSOP. In section \ref{section:4}, the GTSP instances are taken as a benchmark; at first, we test the performance of our meta-heuristic against the optimal solution for the test instances up to 100 nodes with multiple travelers, then we test the meta-heuristic on the medium and large instances up to 500 nodes and 1084 nodes respectively. 217vm1084 is the largest available instance of GTSP. In the mDmSOP, travelers must find the path within a given budget to maximize the collected profit. The results show that if we increase the number of travelers, the profit gained by travelers increases significantly. In section \ref{section:5}, we conclude the paper.

\section{Problem Description and Formulation}
\label{section:2}
The mDmSOP is a generalization of the SOP (Set Orienteering Problem), so first, we give a formal definition of the SOP.

The SOP can be formalized on a directed complete graph $G (V, E)$, where $V=\{\, 1,\ldots,n + m \,\}$ is the set of vertices and $E=\left \{\, (i, j) \mid (i, j) \in V^2 \, \right \}$ is the set of edges. The edge $(i, j)$ is defined as the edge from the vertex $i$ to the vertex $j$. Moreover, a cost $c_{ij}\geq 0$, is associated with the edge $(i, j)$. The vertices are partitioned into disjoint sets $S=\{\, s_{1},\ldots, s_{q} \,\}$ such that their union contains all the vertices of the graph. The objective of the problem is to gain maximum profit by visiting the possible number of sets within a distance constraint ($B$) with the predefined starting and ending node (i.e., node $1$) using a single traveler. The profit from a set can be collected if precisely one node of a set is visited by the traveler. 

In this paper, we propose the mDmSOP, a generalized variant based on the SOP. The mDmSOP can be formally described on a directed complete graph $G (V, E)$, where $V$ is the set of vertices and $E$ is the set of edges, together with the costs $\left \{\, c_{ij} \mid (i, j) \in V^2 \, \right \}$ as defined above. We have $m$ travelers ($t \in \{\, 1,\ldots, m \,\}$) associated with the depots $\{\, n + 1,\ldots,n + m \,\}$. The objective of the mDmSOP is to find $m$ sequences of the clusters using $m$ travelers so that collected profit from all the visited sets can be maximized. The profit from a set can be collected by at most one traveler by visiting one node only. $B$ is the upper distance bound for the sum of the distances traveled by all the travelers, and $B$ is also calculated for the individual traveler in another formulation. The profit of all the depots is $0$.

\begin{figure}[htbp]
    \centerline{\includegraphics[width=0.5\textwidth]{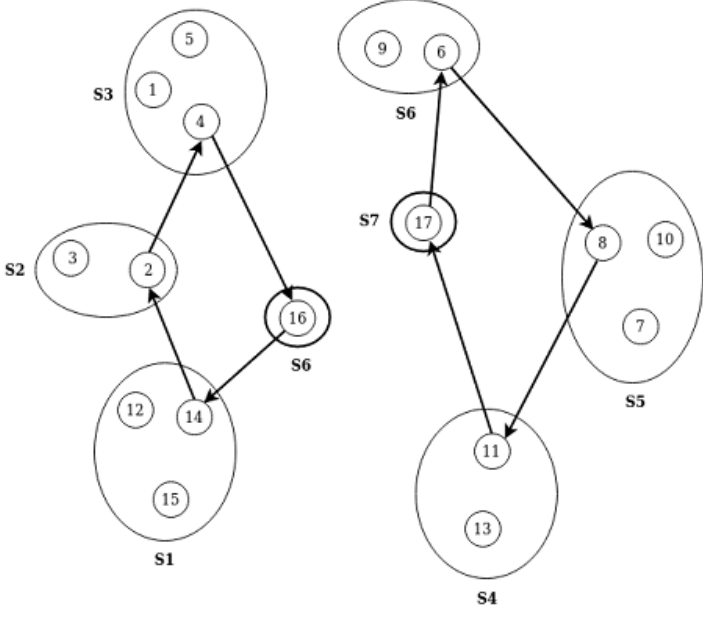}}
    \caption{An example of the mDmSOP.}
    \label{fig:mDmSOP}
\end{figure}

Here, Figure \ref{fig:mDmSOP} illustrates the example of the solution of an instance of the mDmSOP using two travelers. Depots $s_6$ and $s_7$ represent the starting and ending points for traveler $t_1$ and $t_2$, respectively, traveler $t_1$ visits set $s_1$, $s_2$ and $s_3$ while traveler $t_2$ visits $s_6$, $s_5$ and $s_4$. There can be non-visited sets because of the budget constraint.

To represent an Integer Linear Programming (ILP) formulation for the mDmSOP, we use some notations, which are given as follows:
\begin{itemize}
  \item  $t \in \{\, 1,\ldots,m \,\} $ represents a salesmen.
  \item $i, j \in \{\, 1,\ldots,n \,\} \cup \{\, n + 1,\ldots,n + m \,\} $ represent the list of vertices.
  \item $D = \{\, n + 1,\ldots,n + m \,\}$ represent the depots.
  \item $c_{ij}$ represents the weight of the edge $(i, j)$.
  \item $ q \in \{\, 1,\ldots,r \,\} \cup \{\, r + 1,\ldots,r + m \,\}$ represent the indices of the sets (clusters).
  \item $ \{\, s_{r + 1},\ldots,s_{r + m} \,\}$ represent the sets corresponding to the depots.
  \item  $P_q$ represents the profit associated with a set $s_q$.
  \item $B$ represents the budget.
  \item $s$ represents a subset of $V - D$.
\end{itemize}

We can then construct a mathematical formulation using the decision variables:

\begin{itemize}
    \item $x_{tij} = 1$ if the traveler $t$ uses the edge $(i,j) \in E $, and $0$ otherwise.
    \item $y_{ti} = 1$ if the vertex $i$ is visited by the traveler $t$, and $0$ otherwise.
    \item  $z_{tq} = 1$ if any vertex in set $s_q$ is visited by the traveler $t$, and $0$ otherwise.
    \item  $u_{ti}$ : node potentials for the Sub-tour Elimination Constraints (SEC).
\end{itemize}

The proposed mathematical formulation of the mDmSOP:

\begin{equation}
\text{maximize} \: \sum_{t} \sum_{q} P_q z_{tq}, \label{eq:1}
\end{equation}

subject to:

\begin{equation}
x_{tij},y_{ti},z_{tq} \in \{\, 0,1 \,\},\quad \forall t, \forall i, \forall j, \forall q, \label{eq:2}
\end{equation}

\begin{equation}
 \sum_{t} \sum_{i} \sum_{j} x_{tij} c_{ij} \leq B, \label{eq:3}
\end{equation}

\begin{equation}
 y_{t(t+n)} =1, \ \quad \forall t, \label{eq:4}
\end{equation}

\begin{equation}
 \sum_{i \in V-\{\, j \,\}}  x_{tij}=y_{tj}, \ \quad \forall t, \forall j, \label{eq:5}
\end{equation}

\begin{equation}
 \sum_{i \in V-\{\, j \,\}}  x_{tji}=y_{tj}, \ \quad \forall t, \forall j, \label{eq:6}
\end{equation}

\begin{equation}
 \sum_{i \in s_q}  y_{ti}=z_{tq}, \ \quad \forall t, \forall q, \label{eq:7}
\end{equation}

\begin{equation}
 \sum_{t}  z_{tq}\leq 1, \ \quad \forall q, \label{eq:8}
\end{equation}

\begin{equation}
 \sum_{i \in s,  j \in V-s}  x_{tij}\geq 1, \ \quad \forall t, \forall s\subseteq V-D, s\neq \emptyset. \label{eq:9}
\end{equation}

The equation \eqref{eq:9} is used for removing sub-tours in the path, and it can be replaced with the equations \eqref{eq:10} and \eqref{eq:11} because the above equation has exponential ($O \left( 2^n \right)$) number of sub-tour elimination constraints. So it will not be feasible to solve using GAMS and CPLEX, while the equations \eqref{eq:10} and \eqref{eq:11} have only $O \left( n^2 \right)$ sub-tour elimination constraints. Following are the equations replacing equation \eqref{eq:9} for sub-tour elimination:

\begin{equation}
 1 \leq u_{ti} \leq n, \ \quad \forall t, \forall i, \label{eq:10}
\end{equation}

\begin{equation}
 u_{ti}-u_{tj}+1 \leq n(1-x_{tij}), \ \quad \forall t, \forall i, \forall j \:, \: i \neq j. \label{eq:11}
\end{equation}

The objective function \eqref{eq:1} maximizes the collected profits from the sets visited. Constraints \eqref{eq:2} define the domain of the variables $x_{tij}, y_{ti} \text{, and}\, z_{tq}$. Constraint \eqref{eq:3} ensures that budget $B$ is not exceeded. Constraint \eqref{eq:4} ensures that each associated traveler's starting and ending depot is the same. Equations \eqref{eq:5} and \eqref{eq:6} imply that the in-degree is equal to the out-degree of a node except for the starting and ending depots. Constraint \eqref{eq:7} ensures that a set $s_q$ is visited by a traveler $t$ if any node in the set is visited, and at most one node can be visited per set. Constraint \eqref{eq:8} implies that no set can be visited by more than one traveler, while equations \eqref{eq:10} and \eqref{eq:11} are used for sub-tour elimination in the itinerary.  

We have designed the SECs for the mDmSOP based on the TSP formulation proposed by \cite{gavish1978travelling} also. The equations are as follows:

\begin{equation}
0\leq u_{ij}  \leq(n-m+1)\sum_{t=1}^{m}x_{tij}, \: \quad \forall (i,j) \in (V-D)^2, \label{eq:12}
\end{equation}

\begin{equation}
\sum_{j\in V,i\neq j}u_{ij}  -\sum_{j\in V-D,i\neq j}u_{ji} =\sum_{t=1}^{m}y_{ti}, \: \quad \forall i \in V-D. \label{eq:13}
\end{equation}

The performance based on \cite{gavish1978travelling} is recently analyzed by \cite{oncan2009comparative} for the Asymmetric Travelling Salesman Problem (ATSP). 

We also introduce an equation to check the overall difference in profit collected by the travelers if each traveler has an individual budget than the cumulative budget. The equation to satisfy individual budget constraints is given as follows:

\begin{equation}
 \sum_{i} \sum_{j} x_{tij} c_{ij} \leq B, \: \quad \forall t. \label{eq:14}
\end{equation}

In the above mathematical formulation, equations \eqref{eq:1}-\eqref{eq:11} attempt to find out the optimal path with maximization of the profit using permutation of the sets and the vertices which are to be visited in the specific set. In the case of the VNS, the problem is partitioned into several phases: $(i)$ valid initial solution construction for the VNS; $(ii)$ determination of the order of the selected clusters by VNS; and $(iii)$ selecting a node after choosing a set sequence using Dynamic Programming Paradigm (DPP). These phases are described in detail in section \ref{section:3}. 

\section{Variable Neighbourhood Search for the mDmSOP}
\label{section:3}
In this section, we design a meta-heuristic based on the VNS. VNS was developed by \cite{hansen1997variable} to solve a subset of combinatorial optimization problems. 
VNS is chosen to develop for this model after observing its performance to solve the OP by \cite{sevkli2006variable}, the Dubins Orienteering Problem (DOP) by \cite{pvenivcka2017dubins(a)}, and the Orienteering Problem with Neighborhoods (OPN) by \cite{pvenivcka2017dubins(b)}. Recently, it has been used by \cite{pvenivcka2019variable} to solve the SOP. \\
\indent The challenging part of using VNS in our variant is generating a valid initial solution. A simple greedy approach only guarantees a valid solution sometimes because we need every traveler to be used. Therefore, in our case, we first use our greedy algorithm to generate a potential initial solution. If the solution satisfies the problem's constraints, we move on to the main part of the VNS meta-heuristic. But, if the solution does not work, we use the Hungarian algorithm to generate an initial solution (if one exists). This serves two purposes: a. Generating an initial solution and b. Checking if the problem is solvable or not since the Hungarian algorithm gives an assignment with minimal cost. Once we have the solution, we use a similar greedy approach to optimize the Hungarian solution. In the next phase, we employ \texttt{local\_search} and \texttt{shake} functions to optimize the solution generated.

\subsection{Data Structure to Store the Solution}
We need a way to store the potential solutions to our problem in a data structure. The data structure should be decided in a way that optimizes the run time of the different steps of our meta-heuristic. We define a data structure similar to adjacency list $\text{adj}$, where $\text{adj}[i] \hspace{1.0mm} ( i < m$) stores the list of sets that have to be visited by the traveler $i$ in that order. The last index $\text{adj}[m]$ stores the sets that are not visited by any salesman. Given $\text{adj}$, the cost of the tour can be calculated by using dynamic programming. Let's say for salesman $p$, we have the following order of sets $S_p = {s_1, s_2, ... s_{n_p}}$ of size $n_p$ that are visited. We let $\text{dp}_{p}[i][j] \hspace {1.0mm} (0 < i \leq n_p, 0 \leq j < |s_i|)$ store the minimum cost to reach the $j$'th city (let us denote it by $s_{i,j}$) in the set $s_i$.
\[
    \text{dp}_{p}[i][j] = \min_{\forall k \in [0, |s_{i-1}|)}(dp_{p}[i-1][k]+C[s_{i-1,k}][s_{i,j}]).
\]
The sum of minimum elements in the last row for each salesman is the final cost of all the tours, and $\text{dp}_{p}[0][0] = 1$ is the initial condition since every salesman starts from its depot set with $0$ cost, and depot set contains only $1$ city, the depot.

Once the VNS terminates, we have the order in which the sets have to be visited. We use the above DPP to find the minimum cost to get the order in which the cities must be visited by storing cities along with the minimum cost.

\subsection{Initial Solution Construction}
Due to the complexity of the problem, it is challenging to construct an initial solution (not necessarily optimal) that satisfies all the constraints. Therefore, we combine the Hungarian algorithm and a simple greedy approach to generate the initial solution. Hungarian algorithm was introduced by \cite{kuhn1955hungarian} to solve combinatorial problems related to assignment.

\begin{algorithm} [h!]
\begin{algorithmic}[1]
\STATE \textbf{Output:}{Return a valid solution structure}
\STATE pot$\_$sol = generateGreedySolution()
 \IF {validSolution(pot$\_$sol)}
\STATE    \textbf{return} pot$\_$sol\
\ELSE
\STATE   pot$\_$sol = generateHungarianSolution() 
\ENDIF
 \IF {validSolution(pot$\_$sol)}
 \STATE   \textbf{return} optimizeHungarianSolution(pot$\_$sol)
 \ELSE
 \STATE   emptySolution()
 \ENDIF
 \end{algorithmic}
 \caption{Initial solution generation}
 \label{algorithm:1}
\end{algorithm}

Let us start with the description of the functions discussed in Algorithm \ref{algorithm:1} :
\begin{enumerate}
  \item generateGreedySolution: The greedy approach is similar as in \cite{pvenivcka2019variable} to solve the SOP. We iterate over all the sets, and for each set, we find a position for it in the sequence for each salesman. Essentially, we iterate over all the positions and choose the position $\text{pos}$ such that: 
  \[\text{pos} = \arg \min_{j} \frac{L(\text{adj}_{ij})-L(\text{adj})}{p_i},\]
  where $p_i$ is the profit from the newly added set, without violating the maximum cost ($B$) condition.  $L(\text{adj}_{ij})$ is the cost of the configuration after adding set $i$ at position $j$ (from the start) and $L(\text{adj})$ is the cost associated to the old arrangement.

 \item generateHungarianSolution: We model our problem to the assignment problem. For this, we first define the cost matrix $C$. We define cost $C_{ij}$ as the minimum cost among all the cities in the set $j$ visited by the salesman $i$. Now, we know that $m \leq r$; therefore, we can add $r-m$ dummy salesmen to have a square matrix. Once we have modeled the problem into an assignment problem, the Hungarian algorithm gives the minimum cost assignment. Now, the given assignment gives the set that each salesman has to visit so that the total cost of the visit is minimal. If the cost of assignment is greater than $B$, it means that there are no feasible solutions to our model.
  \item optimizeHungarianSolution: This is similar to the greedy approach. The only difference is that, in this case, we start with some non-empty arrangement generated by the Hungarian algorithm and improve it using the same greedy paradigm.
  \item validSolution: This function is responsible for checking if the solution satisfies all the constraints of the problem.
  \item emptySolution: This function is responsible for giving an empty solution structure if no solution exists for the given problem satisfying all the constraints. 
\end{enumerate}

As discussed earlier, the above algorithm returns a valid solution if one exists or an empty structure, which signifies that no solution exists for that instance of the mDmSOP.

\subsection{Variable Neighbourhood Search}
The VNS can be divided into two subparts:
\begin{itemize}
    \item $\texttt{shake}$: This part is used to find the new neighborhood. Each neighborhood has its own $\texttt{shake}$ operation. It is responsible for the exploration phase of the VNS, ensuring that we are not stuck at the local maxima. In our case, we define two $\texttt{shake}$ operators.
    \item $\texttt{local\_search}$: Once we are inside a neighborhood, we want to move towards the local maximum, i.e., the maximum profit corresponding to that particular neighborhood. In our case, we have two different operators to perform the $\texttt{local\_search}$ operation.
\end{itemize}

Together, the two parts are responsible for exploration and exploitation. Now, there can be different ways to perform these operations. It might be the case that for a given neighborhood, we can repeatedly perform the $\texttt{local\_search}$ operation until we can not find further optimizations in the neighborhood. In our meta-heuristic, we perform a $\texttt{shake}$ operation using the first operator, followed by a $\texttt{local\_search}$ operation with the same operator. This is done until no better solution is found using all the operators.

\subsection{The $\texttt{shake}$ Operation}
The VNS has different neighborhood structures where $\texttt{local\_search}$ operations are performed. We have chosen two operators that $\texttt{shake}$ the current solution to find two ($l_{max} = 2$) new neighborhoods to search for the optimal solution. These are:

\begin{enumerate}
    \item Operator $1$ to get to the neighborhood $1$ $(l = 1)$: In this operation, we select a path $S$ (i.e., a subarray) and an index and move the subpart to that position. We randomly decide which side of the index we want to move the path to. While selecting a path, we ensure that the whole path is taken from the same salesman's sequence of visited sets or non-visited sets. Firstly, two indices $l_1$ and $l_2$, $0 \leq l_1, l_2 \leq m$ are randomly selected. From $l_1$, $i_1$, and $j_1$ are chosen randomly such that $0 \leq i_1 \leq j_1 < \text{Length}(\text{adj}[l_1])$, where $\text{Length}(\text{adj}[l_1])$ is the the number of elements in $\text{adj}[l_1]$. Similarly, $0 \leq i_2 < \text{Length}(l_2)$ is chosen randomly. $S = {\text{adj}[l_1][i_1], ..., \text{adj}[l_1][j_1]}$ forms the path.
    \item Operator $2$ to get to the neighborhood $2$ $(l = 2)$: In this operation, we select two paths, $S_1$ and $S_2$, and exchange them. Again, while choosing paths, we ensure that the whole path is taken from the same salesman's sequence of sets or non-visited sets. We also ensure that the two paths are non-intersecting. Again, indices $l_1$ and $l_2$ are chosen randomly. From $l_1$, $i_1$, and $j_1$ are chosen randomly such that $0 \leq i_1 \leq j_1 < \text{Length}(\text{adj}[l_1])$ and $i_2$, $j_2$ such that $0 \leq i_2 \leq j_2 < \text{Length}(\text{adj}[l_2])$, are chosen from $l_2$. The 2 paths formed are $S_1 = {\text{adj}[l_1][i_1], ..., \text{adj}[l_1][j_1]}$, $S_2 = {\text{adj}[l_2][i_2], ..., \text{adj}[l_2][j_2]}$.
\end{enumerate}

In both the $\texttt{shake}$ operations, we always ensure that all the constraints of our problem are satisfied.

\begin{algorithm}[h!]
\begin{algorithmic}[1]
 \STATE \textbf{Input:} $q$ - List of customer sets, $B$ - Maximum allowed Budget
 \STATE \textbf{Output:} $u$ - Solution path within Budget
\STATE $u\leftarrow $ generateInitialSolution $(q, B)$
\WHILE {\textit{termination$\_$condition is false}}
\STATE    $l = 1$\\ 
            \WHILE {$l \leq l_{max}$}
\STATE    curr$\_$neighbourhood = \texttt{shake} ($u$, $l$)
\STATE     pot$\_$sol = \texttt{local\_search}(curr$\_$neighbourhood, $l$)
            \IF{validSolution(pot$\_$sol) \& profit(pot$\_$sol) $>$ profit($u$)} 
\STATE            $u$ = pot$\_$sol
\STATE            l = 1
\ELSE
\STATE           l++
            \ENDIF     
            \ENDWHILE
        \ENDWHILE
\end{algorithmic}
\caption{VNS for the mDmSOP}
\label{algorithm:2}
\end{algorithm}

\subsection{The $\texttt{local\_search}$ Operation}
Local search is the main search operation used to search the current neighborhood. This is responsible for optimizing the potential solution. Again, we define two operators to search over the neighborhood. These are:
\begin{enumerate}
    \item Operator $1$ to search the neighborhood $(l = 1)$: Select two sets $i_1$ and $i_2$ randomly such that\\ $0 \leq i_1, i_2 < n+m$. Now we can put:
    \begin{itemize}
        \item $i_1$ after $i_2$.
        \item $i_1$ before $i_2$.
    \end{itemize}
    This is decided randomly, where each operation can be performed with equal probability.
    \item Operator $2$ to search the neighborhood $(l = 2)$: In this operation, we select two sets randomly and exchange them. Again, $i_1$ and $i_2$ are selected randomly such that $0 \leq i_1, i_2 < n+m$, but in this case, we exchange the sets between these two vertices.
\end{enumerate}

Unlike the $\texttt{shake}$ operation, in the $\texttt{local\_search}$ operation, we do not need to choose the index in $\text{adj}$ since the local search operators are independent of this index. In this case, any two vertices can be selected and replaced. In both cases, we only keep the solution if it increases the value of our profit without violating any of the constraints. Once we have the final solution structure, we use a Dynamic Programming Paradigm (DPP) to find the cost of the solution structure. 

Following are the descriptions of the functions used in Algorithm \ref{algorithm:2}
\begin{enumerate}
\item generateInitialSolution: This function generates an initial solution $u$ within the given budget $B$ using the list of customer sets $q$. Algorithm \ref{algorithm:1} is used to create this initial solution.
\item $\texttt{shake}$: This function generates a neighborhood solution of $u$ by applying a shake procedure. The parameter $l$ controls the intensity of the $\texttt{shake}$ function.
\item $\texttt{local\_search}$: This function performs a local search within the neighborhood defined by the $\texttt{shake}$ procedure. It explores nearby solutions to the current one, trying to improve the solution's quality. 
\item validSolution: This function checks if the potential solution (pot$\_$sol) is valid, meaning it satisfies all the constraints (i.e. budget, etc.).
\item profit: This function calculates the profit associated with the potential solution using DPP by summing up the profits of the visited customer sets.

The algorithm iterates until a termination condition is met, continuously improving the solution $u$ by exploring its neighborhood and performing local searches. If a new solution is found that is both valid and has a higher profit than the current best solution $u$, $u$ is updated to this new solution. Otherwise, the shaking intensity $l$ is increased to explore a broader search space.
\end{enumerate}

\section{Computational Tests}
\label{section:4}
In this section, we analyze the performance of the VNS and the mathematical formulation by giving comparative results for the mDmSOP. The algorithm is implemented in C++ on the Windows 10 platform and runs on an i7-6400 CPU @3.4Ghz processor with 32GB of RAM.

In section \ref{section:4.1}, we describe how the instances are generated for mDmSOP, and the simulation results are shown in section \ref{section:4.2}.

\subsection{Test instances}
\label{section:4.1}
The mDmSOP is a novel variant of the SOP with multiple depots with exactly one traveler associated with each depot. So there is no benchmark instance available for it, which is why we generate some instances adapting the GTSP instances in \cite{noon1988generalized} that are symmetric in nature and follow the Euclidean property. These instances are also used by \cite{archetti2018set} to solve the SOP. We divide our instances into three groups mainly: \\ 1. Set-1: contains small instances of up to $200$ nodes. \\ 2. Set-2: contains mid-size instances of up to $500$ nodes.\\ 3. Set-3: contains large instances of up to $1084$ nodes. \\ 
\indent 217vm1084 is the largest instance available of the GTSP. We performed some changes in the GTSP instances according to the need of our formulation, such as defining the depots. As we used multiple depots and multiple travelers with a single traveler associated with one depot, we started to pick from the last node of the available node in the GTSP instance. We increased it according to our need for depots and made those depots a unique set with one depot only. To generate the profit, we used two methods: 1. The profit generated of a set comes from the function $\text{mod}(C_g)$, where $C_g$ is the number of nodes available in a cluster, and 2. We used a function $(1 + (7141 j))\text{mod}(100)$ used by \cite{archetti2018set} to generate a pseudo-random profit of a node in the set, and then the overall profit associated with a set can be calculated by summing up all the profits of the nodes in that set. These two rules are called $g_1$ and $g_2$, respectively.\\
\indent Budget $B$ can be calculated as $\lceil(w m T_{\text{Max}})\rceil$ for cumulative budget and $\lceil(w T_{\text{Max}})\rceil$ for individual budget, where $T_{\text{Max}}$ is the optimal value for each GTSP instance found by \cite{fischetti1997branch}, 
$w$ is a multiplicative factor so that we can adjust the budget according to our need, and $m$ is the number of travelers used to solve the mDmSOP.\\
\indent So overall, we simulated the results on $20$ small instances in Table \ref{tab: Table 1}. We simulated $51$ small-sized instances in Set $1$, $30$ mid-sized instances in Set $2$, and $27$ large-sized instances up to $1084$ nodes in Set $3$ in Table \ref{tab: Table 2}. In our case, we tested small and mid-size instances with two and three travelers, while for large instances, we took up to four travelers. 

\subsection{Computational results}
\label{section:4.2}
In this section, we show the computational results performed on GTSP instances. We set some termination conditions for the VNS.
\begin{enumerate}
    \item Iteration\_Count $=20$: It corresponds to the number of times an instance is run. 
    \item Max\_Bad\_Iterations$=100000$: It corresponds to the number of iterations in which the algorithm does not improve the solution. 
    \item Max\_Time$=30$ minutes: It corresponds to the maximum time any instance can be run.
\end{enumerate}

Now, we present our results. Table \ref{tab: Table 1} is arranged as the instance name of GTSP, $n$ represents the nodes present in the instance, the number of Travelers, $P_g$ represents the rules to generate the profit, solution and time to solve the instance using CPLEX, solution and time taken by the VNS in the next two columns, and the percentage gap between the mathematical formulation solution given by CPLEX and the meta-heuristic.

We also checked the performance of SECs based on \cite{gavish1978travelling} (equations \eqref{eq:12}-\eqref{eq:13}) in Table \ref{tab: Table 3}, but we observed that except in two out of twelve instances of the GTSP, it takes more time than equations \eqref{eq:10}-\eqref{eq:11} for the mDmSOP.

We divide our experiments into two parts as follows: 
In Table \ref{tab: Table 1}, we compare the mathematical formulation results obtained using CPLEX and the VNS on small instances of the GTSP up to $198$ nodes. The results show that the VNS solves the problem in less time and obtains the optimal solution in all the cases except two instances of 11berlin52 using three travelers and one instance of 16eil76 using two travelers with $g_1$ profit. In these cases, the solutions obtained by our meta-heuristic are $39$, $1630$, and $47$ units with $23.91\%$, $27.68\%$, and $2.08\%$ relative gap, respectively. The CPLEX is not able to solve any instance of 20rat99 and went out of memory (OOM) with the respective time shown in Table \ref{tab: Table 1}. The results also show that as we increased the number of travelers from $2$ to $3$, there is a $38.70\%$ and $44.36\%$ increment in the profit in the case of $g_1$ and $g_2$ while solving the mathematical formulation using CPLEX. In the same scenario, the VNS increased the profit by $32.46\%$ and $35.74\%$ if we generate the profit using rules $g_1$ and $g_2$.

\begin{longtable}[h!]{cccccccccc}
\caption{Comparison of ILP and VNS on the GTSP data instances of the mDmSOP for cumulative budget ($w = 0.25$)} \\
\label{tab: Table 1} \\
\hline \endhead
Instance   & $n$  & Travelers & Pg & \multicolumn{2}{c}{ILP}     &  & \multicolumn{2}{c}{VNS} &            \\ 
\cline{5-6}\cline{8-9}
           &    &            &    & Opt. Solution & Time (sec.) &  & Solution & Time (sec.) & Gap (\%)  \\ 
\hline
11berlin52 & 52 & 2          & $g_1$ & 29            & 850.641     &  & 29       & 5.505       & 0.00      \\
11berlin52 & 52 & 2          & $g_2$ & 1304          & 733.891     &  & 1304     & 5.385       & 0.00      \\
11berlin52 & 52 & 3          & $g_1$ & 46            & 11.469     &  & 39       & 6.350       & 23.91      \\
11berlin52 & 52 & 3          & $g_2$ & 2254          & 226         &  & 1630     & 6.507       & 27.68      \\
11eil51    & 51 & 2          & $g_1$ & 33            & 236.172     &  & 33       & 5.527       & 0.00      \\
11eil51    & 51 & 2          & $g_2$ & 1637          & 330.297     &  & 1637     & 5.502       & 0.00      \\
11eil51    & 51 & 3          & $g_1$ & 41            & 2775.297      &  & 41       & 7.872       & 0.00      \\
11eil51    & 51 & 3          & $g_2$ & 2077          & 700.39       &  & 2077     & 7.032       & 0.00      \\
14st70     & 70 & 2          & $g_1$ & 45            & 147328.391  &  & 45       & 6.962       & 0.00      \\
14st70     & 70 & 2          & $g_2$ & 1966          & 280320.89   &  & 1966     & 6.969       & 0.00      \\
14st70     & 70 & 3          & $g_1$ & 66            & 87439.891   &  & 66       & 11.656       & 0.00      \\
14st70     & 70 & 3          & $g_2$ & 3134          & 36985.969  &  & 3134     & 13.387       & 0.00      \\
16eil76    & 76 & 2          & $g_1$ & 48            & 6966.515    &  & 47       & 7.355       & 2.08      \\
16eil76    & 76 & 2          & $g_2$ & 2324          & 24870.469   &  & 2324     & 7.361       & 0.00      \\
16eil76    & 76 & 3          & $g_1$ & 62            & 206361.641    &  & 62       & 8.924       & 0.00      \\
16eil76    & 76 & 3          & $g_2$ & 2974          & 172772    &  & 2974     & 10.584       & 0.00      \\
20rat99    & 99 & 2          & $g_1$ & OOM         & 38302.985   &  & 46       & 8.158       & -         \\
20rat99    & 99 & 2          & $g_2$ & OOM         & 23520.641   &  & 2077     & 6.702       & -         \\
20rat99    & 99 & 3          & $g_1$ & OOM         & 28064.468   &  & 64       & 10.496       & -         \\
20rat99    & 99 & 3          & $g_2$ & OOM         & 203841.765  &  & 3088     & 9.941       & -         \\
\hline
\end{longtable}

Table \ref{tab: Table 2} is arranged as follows: The GTSP instances are divided into three sets represented by the first column, the next two columns represent the number of nodes present in an instance ($n$), and the number of travelers used to solve the instance. The next four columns represent the solution and time taken by the VNS in the case of profit rule $g_1$ and $g_2$, respectively.
Table \ref{tab: Table 2} is divided as follows:
\begin{enumerate}
    \item Set $1$ contains the GTSP instances of size less than $200$ nodes.
    \item Set $2$ contains the GTSP instances of size from $200$ nodes to $500$ nodes.
    \item Set $3$ contains the GTSP instances of size from $500$ nodes to $1084$ nodes (the largest available GTSP instance).
\end{enumerate}
The average time taken by the meta-heuristic in the case of Set $1$ is $10.567s$ while using two travelers, and the profit earned is $66.038$ units; while using three travelers, the average time is $11.723s$, and the profit is increased by $41.93\%$. A similar pattern is found in the case of rule $g_2$; the profit is increased by $43.17\%$ when using three travelers instead of two. In the case of Set $2$, which takes the GTSP instances of $200$ nodes to $500$ nodes, the average profit earned by the meta-heuristic is $158.133$ units and $245.00$ units while using two and three travelers in case of $g_1$ and $8389.6$ units and $12333.477$ units in case of $g_2$. It is $54.93\%$ and $47.00\%$ profit increment in Set $2$ for rules $g_1$ and $g_2$ when two and three travelers are used, respectively.

We simulate Set $3$ instances for up to four travelers to check the simulation results on large instances for more travelers. The profit increment for rule $g_1$ is $41.21\%$ when we use three travelers instead of two, and $81.35\%$ when we use four travelers instead of three. In the case of rule $g_2$, the profit increases by $28.79\%$ when three travelers are used instead of two, and $97.02\%$ when four travelers are used instead of three.

\begin{longtable}{cccclcclcc}
\caption{Profit comparison on the GTSP data instances of the mDmSOP for cumulative budget ($w=0.25$)} \\
\label{tab: Table 2} \\ 
\hline
              &            &      &            &  & \multicolumn{2}{c}{$g_1$}               &  & \multicolumn{2}{c}{$g_2$}                 \\ 
\cline{6-7}\cline{9-10}
              & Instance   & $n$    & Travelers &  & \multicolumn{2}{c}{VNS}              &  & \multicolumn{2}{c}{VNS}                \\
              &            &      &            &  & Solution          & Time             &  & Solution           & Time              \\
\hline \endhead
Set 1         & 11berlin52 & 52   & 2          &  & 29                & 5.489            &  & 1304               & 5.717             \\
              & 11eil51    & 51   & 2          &  & 33                & 5.622            &  & 1637               & 5.703             \\
              & 14st70     & 70   & 2          &  & 45                & 7.700            &  & 1966               & 6.540             \\
              & 16eil76    & 76   & 2          &  & 48                & 7.585            &  & 2217               & 6.246             \\
              & 20kroA100  & 100  & 2          &  & 49                & 9.165            &  & 2224               & 8.514             \\
              & 20kroB100  & 100  & 2          &  & 60                & 8.476            &  & 3203               & 9.300             \\
              & 20kroC100  & 100  & 2          &  & 48                & 9.494            &  & 2566               & 7.728             \\
              & 20kroD100  & 100  & 2          &  & 63                & 8.507            &  & 3308               & 8.435             \\
              & 20kroE100  & 100  & 2          &  & 59                & 7.635            &  & 3038               & 7.780             \\
              & 20rat99    & 99   & 2          &  & 46                & 7.573            &  & 2077               & 6.383             \\
              & 20rd100    & 100  & 2          &  & 55                & 7.000            &  & 3078               & 8.551             \\
              & 21eil101   & 101  & 2          &  & 75                & 9.234            &  & 3715               & 8.479             \\
              & 21lin105   & 105  & 2          &  & 74                & 8.534            &  & 3673               & 8.786             \\
              & 22pr107    & 107  & 2          &  & 51                & 7.174            &  & 2540               & 7.157             \\
              & 25pr124    & 124  & 2          &  & 65                & 10.635           &  & 3170               & 10.893            \\
              & 26bier127  & 127  & 2          &  & 112               & 15.868           &  & 5713               & 16.467            \\
              & 26ch130    & 130  & 2          &  & 81                & 10.976           &  & 3936               & 10.203            \\
              & 28pr136    & 136  & 2          &  & 61                & 10.127           &  & 3236               & 10.406            \\
              & 29pr144    & 144  & 2          &  & 73                & 12.439           &  & 3622               & 12.847            \\
              & 30ch150    & 150  & 2          &  & 83                & 15.251           &  & 4462               & 13.085            \\
              & 30kroA150  & 150  & 2          &  & 73                & 11.357           &  & 3832               & 11.897            \\
              & 30kroB150  & 150  & 2          &  & 90                & 9.162            &  & 4428               & 12.966            \\
              & 31pr152    & 152  & 2          &  & 76                & 15.875           &  & 3259               & 14.462            \\
              & 32u159     & 159  & 2          &  & 102               & 19.281           &  & 4352               & 16.031            \\
              & 39rat195   & 195  & 2          &  & 85                & 18.521           &  & 4017               & 17.376            \\
              & 40d198     & 198  & 2          &  & 81                & 16.063           &  & 4044               & 19.046            \\
\textbf{Avg.} &            &      &            &  & \textbf{66.038}   & \textbf{10.567}  &  & \textbf{3254.500}  & \textbf{10.423}   \\
              & 11berlin52 & 52   & 3          &  & 35                & 5.724            &  & 1630               & 5.914             \\
              & 11eil51    & 51   & 3          &  & 41                & 6.770            &  & 2077               & 6.480             \\
              & 14st70     & 70   & 3          &  & 66                & 8.976            &  & 3134               & 11.563            \\
              & 16eil76    & 76   & 3          &  & 61                & 9.384            &  & 2974               & 9.166             \\
              & 20kroA100  & 100  & 3          &  & 77                & 8.685            &  & 3806               & 10.019            \\
              & 20kroB100  & 100  & 3          &  & 84                & 9.636            &  & 4242               & 10.474            \\
              & 20kroC100  & 100  & 3          &  & 77                & 10.429           &  & 3883               & 9.267             \\
              & 20kroD100  & 100  & 3          &  & 87                & 8.212            &  & 4321               & 10.631            \\
              & 20kroE100  & 100  & 3          &  & 82                & 7.838            &  & 4042               & 10.035            \\
              & 20rat99    & 99   & 3          &  & 64                & 8.593            &  & 3088               & 10.356            \\
              & 20rd100    & 100  & 3          &  & 83                & 9.638            &  & 4218               & 9.837             \\
              & 21eil101   & 101  & 3          &  & 88                & 9.786            &  & 4411               & 11.460            \\
              & 21lin105   & 105  & 3          &  & 89                & 8.749            &  & 4326               & 9.428             \\
              & 22pr107    & 107  & 3          &  & 50                & 8.187            &  & 2475               & 8.255             \\
              & 25pr124    & 124  & 3          &  & 87                & 16.664           &  & 4367               & 13.698            \\
              & 26bier127  & 127  & 3          &  & 119               & 16.235           &  & 5920               & 16.250            \\
              & 26ch130    & 130  & 3          &  & 114               & 10.028           &  & 5795               & 11.327            \\
              & 28pr136    & 136  & 3          &  & 113               & 15.871           &  & 5661               & 17.286            \\
              & 29pr144    & 144  & 3          &  & 109               & 22.247           &  & 5459               & 18.780            \\
              & 30ch150    & 150  & 3          &  & 126               & 16.975           &  & 6292               & 18.413            \\
              & 30kroA150  & 150  & 3          &  & 109               & 15.023           &  & 5039               & 13.551            \\
              & 31pr152    & 152  & 3          &  & 109               & 9.691            &  & 5395               & 9.664             \\
              & 32u159     & 159  & 3          &  & 134               & 18.965           &  & 6595               & 24.390            \\
              & 39rat195   & 195  & 3          &  & 125               & 13.070           &  & 6238               & 16.449            \\
              & 40d198     & 198  & 3          &  & 181               & 13.993           &  & 9064               & 19.740            \\
\textbf{Avg.} &            &      &            &  & \textbf{93.731}   & \textbf{11.723}  &  & \textbf{4659.654}  & \textbf{12.908}   \\

Set 2         & 40kroa200  & 200  & 2          &  & 104               & 18.175           &  &                  5353               & 19.362            \\
              & 40krob200  & 200  & 2          &  & 105               & 13.841           &  & 5979               & 16.876            \\
              & 45ts225    & 225  & 2          &  & 106               & 15.840           &  & 5542               & 25.063            \\
              & 45tsp225   & 225  & 2          &  & 126               & 29.470           &  & 5988               & 26.103            \\
              & 46pr226    & 226  & 2          &  & 131               & 29.001           &  & 6381               & 28.515            \\
              & 53gil262   & 262  & 2          &  & 146               & 29.165           &  & 7961               & 30.877            \\
              
              & 53pr264    & 264  & 2          &  & 130               & 12.128           &  & 6375               & 12.418            \\
              & 56a280     & 280  & 2          &  & 111               & 22.277           &  & 6041               & 27.823            \\
              & 60pr299    & 299  & 2          &  & 138               & 36.567           &  & 6796               & 32.166            \\
              & 64lin318   & 318  & 2          &  & 193               & 38.068           &  & 9447               & 35.395            \\
              & 80rd400    & 400  & 2          &  & 175               & 50.981           &  & 9639               & 57.933            \\
              & 84fl417    & 417  & 2          &  & 177               & 34.422           &  & 9770               & 61.134            \\
              & 88pr439    & 439  & 2          &  & 230               & 65.001           &  & 12967              & 79.528            \\
              & 89pcb442   & 442  & 2          &  & 137               & 20.129           &  & 8671               & 33.557            \\
              & 99d493     & 493  & 2          &  & 363               & 110.093          &  & 18934              & 147.928           \\
\textbf{Avg.} &            &      &            &  & \textbf{158.133} & \textbf{35.011}  &  & \textbf{8389.6}    & \textbf{42.312}   \\
              & 40kroa200  & 200  & 3          &  & 170               & 26.690           &  & 8392               & 32.456            \\
              & 40krob200  & 200  & 3          &  & 162               & 28.054           &  & 8403               & 33.594            \\
              & 45ts225    & 225  & 3          &  & 170               & 35.485           &  & 8675               & 37.941            \\
              & 45tsp225   & 225  & 3          &  & 171               & 34.553           &  & 9398               & 38.721            \\
              & 46pr226    & 226  & 3          &  & 189               & 40.492           &  & 9361               & 39.425            \\
              & 53gil262   & 262  & 3          &  & 218               & 44.732           &  & 10959              & 65.137            \\
              & 53pr264    & 264  & 3          &  & 175               & 31.049           &  & 9598               & 53.367            \\
              & 56a280     & 280  & 3          &  & 198               & 22.514           &  & 10361              & 30.629            \\
              & 60pr299    & 299  & 3          &  & 240               & 52.445           &  & 12152              & 49.874            \\
              & 64lin318   & 318  & 3          &  & 261               & 71.079           &  & 12382              & 60.482            \\
              & 80rd400    & 400  & 3          &  & 291               & 105.174          &  & 15517              & 98.928            \\
              & 84fl417    & 417  & 3          &  & 324               & 124.175          &  & 15666              & 135.842           \\
              & 88pr439    & 439  & 3          &  & 364               & 123.366          &  & 18116              & 92.200            \\
              & 89pcb442   & 442  & 3          &  & 279               & 44.007           &  & 13026              & 65.144            \\
              & 99d493     & 493  & 3          &  & 463               & 163.741          &  & 22996              & 160.352           \\
\textbf{Avg.} &            &      &            &  & \textbf{245}      & \textbf{63.170}  &  & \textbf{12333.467} & \textbf{66.273}   \\

Set 3         & 115rat575  & 575  & 2          &  & 213               & 26.990           &  & 11163               & 29.747            \\
              & 115u574    & 574  & 2          &  & 232               & 60.712           &  & 11704              & 58.485            \\
              & 131p654    & 654  & 2          &  & 301               & 101.929          &  & 15283              & 106.203           \\
              & 132d657    & 657  & 2          &  & 327               & 109.827          &  & 14464              & 83.309            \\
              & 145u724    & 724  & 2          &  & 201               & 37.576           &  & 11234               & 43.747            \\
              & 157rat783  & 783  & 2          &  & 237               & 30.534           &  & 11785              & 26.547            \\
              & 201pr1002  & 1002 & 2          &  & 324               & 95.345           &  & 15275              & 136.370           \\
              & 212u1060   & 1060 & 2          &  & 312               & 120.284          &  & 14182              & 103.254           \\
              & 217vm1084  & 1084 & 2          &  & 352               & 105.873          &  & 22907              & 177.420           \\
\textbf{Avg.} &            &      &            &  & \textbf{277.667}  & \textbf{76.563}  &  & \textbf{14221.889} & \textbf{85.009}   \\
  
              & 115rat575  & 575  & 3          &  & 328               & 64.201           &  & 17369              & 84.721            \\
              & 115u574    & 574  & 3          &  & 393               & 112.916          &  & 17075              & 94.751            \\
              & 131p654    & 654  & 3          &  & 426               & 199.741          &  & 20784              & 163.031           \\
              & 132d657    & 657  & 3          &  & 369               & 122.919          &  & 17774              & 85.449            \\
              & 145u724    & 724  & 3          &  & 316               & 64.340           &  & 13501              & 45.544            \\
              & 157rat783  & 783  & 3          &  & 432               & 67.672           &  & 20365              & 79.029            \\
              & 201pr1002  & 1002 & 3          &  & 423               & 141.951          &  & 18931              & 122.821           \\
              & 212u1060   & 1060 & 3          &  & 306               & 36.601           &  & 14792              & 32.046            \\
              & 217vm1084  & 1084 & 3          &  & 536               & 231.019          &  & 24259              & 207.371           \\
\textbf{Avg.} &            &      &            &  & \textbf{392.111}      & \textbf{115.706} &  & \textbf{18316.667} & \textbf{101.640}   \\
              & 115rat575  & 575  & 4          &  & 533               & 133.318          &  & 26938              & 142.053           \\
              & 115u574    & 574  & 4          &  & 550               & 126.553          &  & 27178              & 140.103           \\
              & 131p654    & 654  & 4          &  & 641               & 143.053          &  & 31852              & 159.132           \\
              & 132d657    & 657  & 4          &  & 649               & 144.902          &  & 32188              & 147.204           \\
              & 145u724    & 724  & 4          &  & 694               & 147.931          &  & 34342              & 161.872           \\
              & 157rat783  & 783  & 4          &  & 724               & 210.805          &  & 36179              & 153.247           \\
              & 201pr1002  & 1002 & 4          &  & 776               & 190.231          &  & 41299              & 264.575           \\
              & 212u1060   & 1060 & 4          &  & 965               & 294.131          &  & 48444              & 421.228           \\
              & 217vm1084  & 1084 & 4          &  & 868               & 367.208          &  & 46373              & 499.003           \\
\textbf{Avg.} &            &      &            &  & \textbf{711.111}  & \textbf{195.348} &  & \textbf{36088.111} & \textbf{232.046} \\
\hline
\end{longtable}

\begin{longtable}[h!]{cccccc}
\caption{Optimal Solution and Time of ILP (Gavish-based) for cumulative budget ($w=0.25$)} \\
\label{tab: Table 3} \\
\hline
Instance   & $n$  & Travelers & Pg & \multicolumn{2}{c}{ILP}     \\ \cline{5-6} 
           &    &            &    & Opt. Solution & Time (sec.) \\ \hline \endhead
11berlin52 & 52 & 2          & $g_1$ & 29            & 161.437     \\
11berlin52 & 52 & 3          & $g_1$ & 46            & 287.375     \\
11berlin52 & 52 & 2          & $g_2$ & 1304          & 541.672     \\
11berlin52 & 52 & 3          & $g_2$ & 2254          & 2352.187    \\
11eil51    & 51 & 2          & $g_1$ & 33            & 5615.985    \\
11eil51    & 51 & 3          & $g_1$ & 41            & 8669.75     \\
11eil51    & 51 & 2          & $g_2$ & 1637          & 6968.672    \\
11eil51    & 51 & 3          & $g_2$ & 2077          & 12550.875   \\
16eil76    & 76 & 2          & $g_1$ & OOM         & 580430.969  \\
16eil76    & 76 & 3          & $g_1$ & 62            & 1311150.891 \\
16eil76    & 76 & 2          & $g_2$ & 2324          & 677789.391  \\
16eil76    & 76 & 3          & $g_2$ & 2974          & 214940.719  \\ \hline
\end{longtable}

\begin{longtable}[h!]{cccccccccc}
\caption{Comparison of ILP and VNS on the GTSP data instances of the mDmSOP for individual budget ($w = 0.25$)} \\
\label{tab: Table 4} \\ 
\hline
Instance   & $n$  & Travelers & Pg & \multicolumn{2}{c}{ILP}     &  & \multicolumn{2}{c}{VNS}            \\ 
\cline{5-6}\cline{8-9}
           &    &            &    & Opt. Solution & Time (sec.) &  & Solution & Time (sec.) & Gap (\%)  \\ 
\hline \endhead
11berlin52 & 52 & 2          & $g_1$ & 27            & 12.75       &  & 27       & 6.296       & 0.00      \\
11berlin52 & 52 & 2          & $g_2$ & 1276          & 12.532      &  & 1276     & 6.078       & 0.00      \\
11berlin52 & 52 & 3          & $g_1$ & 38            & 131.328     &  & 38       & 6.709       & 0.00      \\
11berlin52 & 52 & 3          & $g_2$ & 1816          & 70.688      &  & 1816     & 7.190       & 0.00      \\
11eil51    & 51 & 2          & $g_1$ & 28            & 6.812       &  & 28       & 6.051       & 0.00      \\
11eil51    & 51 & 2          & $g_2$ & 1552          & 7.625       &  & 1552     & 5.929       & 0.00      \\
11eil51    & 51 & 3          & $g_1$ & 37            & 54.063      &  & 37       & 5.974       & 0.00      \\
11eil51    & 51 & 3          & $g_2$ & 1862          & 52.11       &  & 1862     & 6.154       & 0.00      \\
14st70     & 70 & 2          & $g_1$ & 26            & 42362.5	   &  & 26       & 6.367       & 0.00      \\
14st70     & 70 & 2          & $g_2$ &  1246         & 56815.047    &  & 1246    & 6.498       & 0.00      \\
14st70     & 70 & 3          & $g_1$ &  OOM        & 428630.157   &  & 47      & 7.825       & -      \\
14st70     & 70 & 3          & $g_2$ &  OOM        & 73645.922    &  & 2398	 & 8.061       & -      \\
16eil76    & 76 & 2          & $g_1$ & 44            & 39798.297    &  & 44      & 6.373       & 0.00      \\
16eil76    & 76 & 2          & $g_2$ & 2136          & 41993.75	    &  & 2136    & 6.404       & 0.00      \\
16eil76    & 76 & 3          & $g_1$ & OOM         & 412576.109   &  & 54      & 6.689       & -      \\
16eil76    & 76 & 3          & $g_2$ & OOM         & 351587.687   &  & 2693    & 6.959       & -      \\
20rat99    & 99 & 2          & $g_1$ & OOM         & 64333.828    &  & 27      & 6.554       & -         \\
20rat99    & 99 & 2          & $g_2$ & OOM         & 28214.422   &  & 1323     & 6.841       & -         \\
20rat99    & 99 & 3          & $g_1$ & OOM         & 37397.156   &  & 32       & 7.301       & -         \\
20rat99    & 99 & 3          & $g_2$ & OOM         & 28017.891   &  & 1643     & 7.571       & -         \\
\hline
\end{longtable}

\begin{longtable}[h!]{cccclcclcc}
\caption{Profit comparison on the GTSP data instances of the mDmSOP for individual budget ($w=0.25$)} \\
\label{tab: Table 5} \\ 
\hline
              &            &      &            &  & \multicolumn{2}{c}{$g_1$}               &  & \multicolumn{2}{c}{$g_2$}                 \\ 
\cline{6-7}\cline{9-10}
              & Instance   & $n$    & Travelers &  & \multicolumn{2}{c}{VNS}              &  & \multicolumn{2}{c}{VNS}                \\
              &            &      &            &  & Solution          & Time             &  & Solution           & Time              \\ 
\hline \endhead

Set 1         & 11berlin52 & 52   & 2          &  & 27                & 6.296            &  & 1276               & 6.078             \\
              & 11eil51    & 51   & 2          &  & 28                & 6.051            &  & 1552               & 5.929             \\
              & 14st70     & 70   & 2          &  & 26                & 6.367            &  & 1246               & 6.498             \\
              & 16eil76    & 76   & 2          &  & 44                & 6.373            &  & 2136               & 6.404             \\
              & 20kroA100  & 100  & 2          &  & 41                & 6.220            &  & 2032               & 6.579             \\
              & 20kroB100  & 100  & 2          &  & 48                & 7.082            &  & 2632               & 7.669             \\
              & 20kroC100  & 100  & 2          &  & 41                & 6.431	           & 
              & 2022	            &  6.911             \\
              & 20kroD100  & 100  & 2          &  & 48	            & 6.811	           & 
              & 2261	            & 7.347              \\
              & 20kroE100  & 100  & 2          &  & 54	            & 7.832	           &
              & 2704	            & 7.969              \\
              	
              & 20rat99    & 99   & 2          &  & 27	            & 6.554            & 
              & 1323	            & 6.841               \\
              & 20rd100    & 100  & 2          &  & 48	& 6.397	     & & 2444	& 6.517                          \\
              & 21eil101   & 101  & 2          &  & 71	& 7.839      & & 3513	& 8.585                          \\
              & 21lin105   & 105  & 2          &  & 73	& 9.949	     & & 3673	& 10.596                          \\
              & 22pr107    & 107  & 2          &  & 21	& 10.070     & & 1035	& 9.903                          \\
              & 25pr124    & 124  & 2          &  & 39	& 8.428	     & & 2029	& 8.160                          \\
              & 26bier127  & 127  & 2          &  & 108	& 14.516	& & 5562	& 13.287                          \\
              & 26ch130    & 130  & 2          &  & 74	& 11.176	& & 3655	& 11.349                  \\
              & 28pr136    & 136  & 2          &  & 44	& 8.016	    & & 2199	& 7.519                  \\
              & 29pr144    & 144  & 2          &  & 42	& 7.442	    & &2067	    & 8.230                 \\
              & 30ch150    & 150  & 2          &  & 70	& 10.316	& & 3350	& 9.576                  \\
              & 30kroA150  & 150  & 2          &  & 60	& 8.160 	& & 2969	& 7.971                  \\
              & 30kroB150  & 150  & 2          &  & 87	& 11.470	& & 4417	& 13.483                  \\
              & 31pr152    & 152  & 2          &  & 35	& 6.338	   & & 1711	    & 6.541                   \\
              & 32u159     & 159  & 2          &  & 76	& 11.556	& &3842	&13.335                  \\
              & 39rat195   & 195  & 2          &  & 54	& 9.049	    & &2658	&9.010                  \\
              & 40d198     & 198  & 2          &  & 24	& 13.465	& &1288	&13.346                   \\
\textbf{Avg.} &            &      &            &  & \textbf{50.385}   & \textbf{8.469}  &  & \textbf{2522.923}  & \textbf{8.678}   \\

              & 11berlin52 & 52   & 3          &  & 38	&6.709	& &1816	&7.190  \\
              & 11eil51    & 51   & 3          &  & 37	&5.974	& &1862	&6.154 \\
              & 14st70     & 70   & 3          &  & 47	&7.825	& &2398	&8.061 \\
              & 16eil76    & 76   & 3          &  & 54	&6.689	& &2693	&6.959\\
              & 20kroA100  & 100  & 3          &  & 69	&8.072	& &3469	&9.505 \\
              & 20kroB100  & 100  & 3          &  & 57	&7.976	& &3036	&8.316 \\
              & 20kroC100  & 100  & 3          &  & 60	&8.707	& &2882	&9.495\\
              & 20kroD100  & 100  & 3          &  & 69	&7.755	& &3213	&8.098\\
              & 20kroE100  & 100  & 3          &  & 72	&10.007	& &3721	&8.582\\
              & 20rat99    & 99   & 3          &  & 32	&7.301	& &1643	&7.571 \\
              & 20rd100    & 100  & 3          &  & 67	&7.615	& &3387	&7.933   \\
              & 21eil101   & 101  & 3          &  & 81	&8.141	& &4033	&8.295   \\
              & 21lin105   & 105  & 3          &  & 80	&8.850	& &3912	&9.138   \\
              & 22pr107    & 107  & 3          &  & 26	&7.553	& &1279	&7.863   \\
              & 25pr124    & 124  & 3          &  & 38	&12.196	& &1967	&12.660   \\
              & 26bier127  & 127  & 3          &  & 118	&15.324	& &5877	&15.863   \\
              & 26ch130    & 130  & 3          &  & 84	&11.452	& &4261	&11.258 \\
              & 28pr136    & 136  & 3          &  & 75	&10.972	& &3643	&9.895 \\
              & 29pr144    & 144  & 3          &  & 47	&7.446	& &2266	&7.788   \\
              & 30ch150    & 150  & 3          &  & 96	&14.046	& &4640	&14.168   \\
              & 30kroA150  & 150  & 3          &  & 87	&13.202	& &4238	&14.054   \\
              & 30kroB150  & 150  & 3          &  & 109	&12.328	& &5483	&14.297 \\
              & 31pr152    & 152  & 3          &  & 34	&8.573	& &1701	&8.830   \\
              & 32u159     & 159  & 3          &  & 87	&11.869	& &4268	&12.468   \\
              & 39rat195   & 195  & 3          &  & 64	&7.516	& &3105	&7.733   \\
              & 40d198     & 198  & 3          &  & 23	&15.205	& &1192	&15.617   \\	
\textbf{Avg.} &            &      &            &  & \textbf{63.500}   & \textbf{9.589}  &  & \textbf{3153.269}  & \textbf{9.915}   \\

Set 2         & 40kroa200  & 200  & 2          &  & 95	&15.202	&&5044	&16.721            \\
              & 40krob200  & 200  & 2          &  & 109	&14.337	&&5271	&14.949          \\
              & 45ts225    & 225  & 2          &  & 87	&14.663	&&4321	&16.012          \\
              & 45tsp225   & 225  & 2          &  & 104	&15.066	&&5218	&15.982          \\
              
              & 46pr226    & 226  & 2          &  & 75	&8.758	&&3707	&8.718   \\
              & 53gil262   & 262  & 2          &  & 136	&21.396	&&6817	&21.710  \\
              
              & 53pr264    & 264  & 2          &  & 130	&11.501	&&6375	&11.402   \\
              & 56a280     & 280  & 2          &  & 108	&10.163	&&5166	&9.904  \\
              & 60pr299    & 299  & 2          &  & 135	&25.517	&&6951	&25.206   \\
              
        & 64lin318   & 318  & 2          &  & 175	&41.482	&&8922	&34.431   \\
              & 80rd400    & 400  & 2          &  & 169	&33.708	&&9143	&39.976   \\
              & 84fl417    & 417  & 2          &  & 195	&44.810	&&9787	&41.420   \\
              & 88pr439    & 439  & 2          &  & 80	&22.888	&&4200	&27.297   \\
              & 89pcb442   & 442  & 2          &  & 160	&30.545	&&8013	&35.566   \\
              & 99d493     & 493  & 2          &  & 370	&113.273	&&18425	 &116.290   \\	
\textbf{Avg.} &            &      &            &  & \textbf{141.867} & \textbf{28.220}  &  & \textbf{7157.333}    & \textbf{29.039}   \\
              & 40kroa200  & 200  & 3          &  & 136	&14.274	&&6899	&13.887   \\
              & 40krob200  & 200  & 3          &  & 143	&21.830	&&7257	&18.150   \\
              & 45ts225    & 225  & 3          &  & 98	&11.565	&&4855	&11.788   \\
              & 45tsp225   & 225  & 3          &  & 144	&21.494	&&6951	&21.499   \\
              & 46pr226    & 226  & 3          &  & 74	&10.411	&&3663	&10.380   \\
              & 53gil262   & 262  & 3          &  & 184	&17.988	&&9937	&25.533   \\
              & 53pr264    & 264  & 3          &  & 129	&11.612	&&6373	&11.535   \\
              & 56a280     & 280  & 3          &  & 125	&10.452	&&6041	&10.427   \\
              & 60pr299    & 299  & 3          &  & 155	&33.690	&&7776	&36.135   \\
              & 64lin318   & 318  & 3          &  & 228	&52.775	&&12007	&62.262   \\
              & 80rd400    & 400  & 3          &  & 275	&67.514	&&13091	&56.916   \\
              & 84fl417    & 417  & 3          &  & 198	&46.483	&&9851	&44.808   \\
              & 88pr439    & 439  & 3          &  & 114	&21.217	&&4408	&19.274   \\
              & 89pcb442   & 442  & 3          &  & 247	&62.587	&&12930	&62.321   \\
              & 99d493     & 493  & 3          &  & 403	&125.441	&&19848	&106.099 \\	
\textbf{Avg.} &            &      &            &  & \textbf{176.867}      & \textbf{35.289}  &  & \textbf{8792.467} & \textbf{34.068}   \\

Set 3         & 115rat575  & 575  & 2          &  & 190	  &18.220	&&9978	&19.258 \\
              & 115u574    & 574  & 2          &  & 192	&33.994	  &&9276  &34.037   \\
              & 131p654    & 654  & 2          &  & 244	&114.935 &&15232	&120.782   \\
              & 132d657    & 657  & 2          &  & 268	&70.209	&&13739	&96.086   \\
              & 145u724    & 724  & 2          &  & 219	&32.060	&&11427	&48.600   \\
              & 157rat783  & 783  & 2          &  & 225	&23.167	&&11678	&25.809   \\
              & 201pr1002  & 1002 & 2          &  & 226	&74.501	&&13940	&58.274   \\
              & 212u1060   & 1060 & 2          &  & 296	&87.102	&&9148	&71.830   \\
              & 217vm1084  & 1084 & 2          &  & 386	&162.132	&&20829	&184.515   \\
\textbf{Avg.} &            &      &            &  & \textbf{249.556	}  & \textbf{68.480}  &  & \textbf{12805.222} & \textbf{73.243}   \\

              & 115rat575  & 575  & 3          &  & 260	&20.845	    &&13606	&26.526   \\
              & 115u574    & 574  & 3          &  & 251	&37.797	    &&12710	&39.724   \\
              & 131p654    & 654  & 3          &  & 315	&116.110	&&15493	&119.954   \\
              & 132d657    & 657  & 3          &  & 353	&68.105	    &&17509	&62.681   \\
              & 145u724    & 724  & 3          &  & 300	&67.783	    &&15926	&73.857   \\
              & 157rat783  & 783  & 3          &  & 371	&38.071	    &&17659	&32.664   \\
              & 201pr1002  & 1002 & 3          &  & 293	&82.344 	&&15637	&110.125   \\
              & 212u1060   & 1060 & 3          &  & 319	&64.135	    &&15000	&49.483   \\
              & 217vm1084  & 1084 & 3          &  & 578	&312.781	&&30881	&307.218 \\	
\textbf{Avg.} &            &      &            &  & \textbf{337.778}      & \textbf{89.775} &  & \textbf{17157.889} & \textbf{91.359}   \\
              & 115rat575  & 575  & 4          &  & 310	&23.291	&&15523	&22.613   \\
              & 115u574    & 574  & 4          &  & 290	&24.381	&&14099	&26.043  \\
              & 131p654    & 654  & 4          &  & 315	&109.630 &&15728	&128.146   \\
              & 132d657    & 657  & 4          &  & 430	&76.379	&&20693	&78.726   \\
              & 145u724    & 724  & 4          &  & 394	&46.047	&&18709	&55.462   \\
              & 157rat783  & 783  & 4          &  & 451	&31.817	&&22831	&38.682   \\
              & 201pr1002  & 1002 & 4          &  & 695	&179.286	&&30469	&166.337   \\
              & 212u1060   & 1060 & 4          &  & 427	&50.454	    &&22059	&62.098   \\
              & 217vm1084  & 1084 & 4          &  & 759	&508.563	&&35258	&387.317 \\	
\textbf{Avg.} &            &      &            &  & \textbf{452.333}  & \textbf{116.650} &  & \textbf{21707.667} & \textbf{107.269} \\
\hline
\end{longtable}

The structure of Table \ref{tab: Table 4} and Table \ref{tab: Table 5} is the same as Table \ref{tab: Table 1} and Table \ref{tab: Table 2}, respectively. These tables show the results of the mDmSOP when each traveler has an individual budget constraint. The results of Table \ref{tab: Table 4} represent that it is evident that the VNS performs far better than the ILP to find out the optimal solution for all the twelve instances which CPLEX is able to solve for individual budget constraints, while CPLEX is not able to solve the other eight instances optimally and goes OOM. The results show that in the case of cumulative budget constraint, the ILP and the VNS are able to gain more profit from clusters because it is more flexible to share the budget among travelers. The overall profit of Set $1$ for the cumulative budget constraint is $31.07\%$ and $28.99\%$ more than that of the individual budget constraint for the mDmSOP for $g_1$ and $g_2$, respectively, using two travelers. In the same situation, if we use three travelers, the relative profit for cumulative budget constraint increases to $47.60\%$ and $47.77\%$. For Set $2$, it is $11.47\%$ \& $17.22\%$ and $38.52\%$ \& $40.27\%$ more than individual budget constraint for two and three travelers applying rule $g_1$ and $g_2$. In the case of Set $3$ instances, the relative profit between cumulative budget constraint and individual budget constraint is $11.26\%$ \& $11.06\%$, $16.09\%$ \& $6.75\%$, and $57.20\%$ \& $66.25\%$ for $g_1$ and $g_2$ respectively, while using two, three and four travelers. 

\section{Conclusion}
\label{section:5}
In this paper, we introduced a new variant of the Set Orienteering Problem, which has multiple depots and the same numbers of travelers, with exactly one traveler associated with each depot; the objective of this problem is to gain maximum profit out of mutually exclusive sets by using multiple travelers with fixed starting and ending points within a fixed budget $B$. The budget $B$ is taken as cumulative as well as individual (i.e. for each traveler). A fixed profit is associated with each cluster, calculated using two rules named $g_1$ and $g_2$, and the profit can only be gained if any traveler visits exactly one node of a set. The results show that the meta-heuristic produces high-quality results as compared to CPLEX, which has less computational time. 

This problem is an extension of the multi-depot multiple traveling salesman problem. It has an application in the supply chain, where a distributor has more than one service point from which the distributor can supply the products to the retailers and gain the maximum profit while giving them a better product price within a given time.

\bibliographystyle{apalike}
\bibliography{main}

\end{document}